\documentclass[12pt,A4]{article}
\usepackage{theorem}
\newtheorem{lemma}{Lemma}
\newtheorem{proposition}[lemma]{Proposition}
\newtheorem{theorem}[lemma]{Theorem}
\newtheorem{corollary}[lemma]{Corollary}
{\theorembodyfont{\upshape}}
{\theorembodyfont{\upshape}}
{\theorembodyfont{\upshape}\newtheorem{note}[lemma]{Note}}
{\theorembodyfont{\upshape}\newtheorem{example}[lemma]{Example}}
{\theorembodyfont{\upshape}}
\newcommand{\Z}{{\bf Z}}
\newcommand{\R}{{\bf R}}
\newcommand{\C}{{\bf C}}
\newcommand{\rme}{{\rm e}}

\newcommand{\cC}{{\cal C}}
\newcommand{\cD}{{\cal D}}

\newcommand{\cH}{{\cal H}}

\newcommand{\cP}{{\cal P}}

\newcommand{\sig}{\sigma}
\newcommand{\alp}{\alpha}
\newcommand{\bet}{\beta}

\newcommand{\kap}{\kappa}
\newcommand{\lam}{\lambda}
\newcommand{\del}{\delta}
\newcommand{\eps}{\varepsilon}
\newcommand{\ome}{\omega}

\newcommand{\Ome}{\Omega}
\newcommand{\Spec}{{\rm Spec}}

\newcommand{\norm}{\Vert}

\newcommand{\Proof}{\underbar{Proof}{\hskip 0.1in}}

\newcommand{\Num}{{\rm Num}}
\newcommand{\conv}{{\rm conv}}
\newcommand{\dist}{{\rm dist}}
\newcommand{\Schrodinger}{Schr\"odinger }
\newcommand{\imp}{\Rightarrow}

\newcommand{\sigl}{\sig_{{\rm loc}}}
\newcommand{\var}{{\rm var}}
\newcommand{\tensor}{\otimes}
\title{SPECTRAL PROPERTIES OF\\
RANDOM NON-SELF-ADJOINT\\
 MATRICES AND OPERATORS}
\author{E.B. Davies}
\date{February 2000}
\begin{document}
\maketitle
%%%%%%%%%%%%%%%%%%%%%
\begin{abstract}
We describe some numerical experiments which determine the degree
of spectral instability of medium size randomly generated matrices
which are far from self-adjoint. The conclusion is that the
eigenvalues are likely to be intrinsically uncomputable for similar
matrices of a larger size. We also describe a stochastic family of
bounded operators in infinite dimensions for almost all of which
the eigenvectors generate a dense linear subspace, but the
eigenvalues do not determine the spectrum. Our results imply that
the spectrum of the non-self-adjoint Anderson model changes
suddenly as one passes to the infinite volume limit.
\vskip 0.1in
AMS subject classifications: 65F15, 65F22, 15A18, 15A52, 47A10,
47A75, 47B80, 60H25.
\par
keywords: eigenvalues, spectral instability, matrices,
computability, pseudospectrum, \Schrodinger operator, Anderson
model.
\end{abstract}
%%%%%%%%%%%%%%%%%%%%%%%%
\section{Introduction}
\par
In a series of recent papers [2,4-10] a number of authors have
investigated the spectral properties of non-self-adjoint matrices
and operators, coming to the conclusion that the eigenvalues are
frequently highly unstable under small perturbations of the
coefficients of the matrices. Trefethen has investigated a series
of numerical examples using the concept of pseudospectrum (the
contour plots of the resolvent norm), which provides a graphical
demonstration of the degree of instability, \cite{Tre1,Tre2}. This
concept is, however, not well adapted to the consideration of very
large numbers of randomly generated matrices, for which one needs
to produce a numerical measure of instability of the spectrum. In
this paper we define such an instability index, and compute it for
a series of randomly generated $N\times N$ matrices for various
values of $N$ up to $50$. Our numerical results are presented in
Section 3, following a short theoretical section which describes
the concepts involved. Our results are fully in line with what
would be expected by experts in pseudospectral theory, but we
believe that such a systematic quantitative investigation of
non-self-adjoint random matrices has not previously been carried
out, and it is clear that much of the spectral theory community is
not aware of these phenomena.

From Section 5 onwards we consider a related problem for a
stochastic family of non-self-adjoint bounded operators in infinite
volume, i.e acting on $l^2(\Z^n)$. For one example we prove that
although the eigenvectors of almost all of the operators span a
dense linear subspace, they do not form a basis, and the spectrum
is much larger than the closure of the set of eigenvalues. Section
6 is devoted to spelling out the implications of our results for
the non-self-adjoint Anderson model of Hatano-Nelson, which has
been the focus of much recent attention.\cite{Gold,
Hatano1,Hatano2,Nelson} We find that the asymptotic behaviour of
the eigenvalues as the volume increases does not describe the full
spectrum of the infinite volume problem. The reason is that there
are many approximate eigenvalues of the finite volume problem which
are not close to true eigenvalues but which nevertheless affect the
infinite volume limit. This situation is typical of
non-self-adjoint operators, for which the eigenvectors need not be
even approximately orthogonal.
\par

\section{The Theoretical Context}

Throughout this section we suppose that $A$ is an $N\times N$
matrix with distinct eigenvalues. This is generically true, since
the set of matrices with repeated eigenvalues forms a lower
dimensional set with zero Lebesgue measure. All vectors below are
assumed to be column vectors unless otherwise stated, and ${}^\ast$
denotes the conjugate transpose. Norms of vectors are always their
Euclidean norms and norms of operators are always the corresponding
operator norms.

Let $\lam$ be an eigenvalue of $A$ of multiplicity $1$, with
corresponding eigenvector $\phi$. Let $\psi$ be an eigenvector of
$A^\ast$ with eigenvalue $\mu$. If $\mu\not= \overline{\lam}$ then
$\psi^\ast\phi=0$ and otherwise this inner product is non-zero.
Assuming that $\mu=\overline{\lam}$, the spectral projection $P$
associated to the eigenvalue is given by
\[
P f=(\psi^\ast\phi)^{-1}(\psi^\ast f)\phi
\]
and its norm is
\begin{equation}
\norm P \norm =\norm \phi \norm \, \norm \psi
\norm/|\psi^\ast\phi |.\label{Pnorm}
\end{equation}

If we label all the above quantities by $n$ for $n=1,..., N$, then
we may define the instability index of $A$ by
\[
i(A)=\max\{ \norm P_n \norm:n=1,...,N\}.
\]
This is unrelated to its condition number
\[
\kap (A)=\norm A \norm \, \norm A^{-1} \norm
\]
since the instability index of a self-adjoint or normal matrix
always equals $1$, but its condition number may be arbitrarily
large; in the converse direction the matrix
\[
\left( \begin{array}{cc}
1&1\\ 0&1+\del
\end{array}
\right)
\]
has small condition number but unbounded instability index as
$\del\to 0$. The norm of a particular spectral projection is also
called the condition number of the eigenvalue, and is known to
measure how unstable the eigenvalue is under small perturbations of
the matrix, \cite{AD,Wat,Wil},\cite[sect. 11.2]{LT}. We emphasize
that if the norm of the spectral projection is very large the
instability of the eigenvalue it intrinsic: it does not depend on
the particular method of computing it. The norm of $P$ is always at
least $1$ and equals $1$ if and only if $P$ is orthogonal, or
equivalently if $\phi=\psi$. The following proposition relates this
index to other measures of how far the matrix is from being normal,
\cite{AD,Wat,Wil}.

\begin{proposition}
We have the relations
\[
(i)\imp (ii) \imp (iii) \imp (iv)
\]
between the following conditions, the constant $k$ being the same
in all cases.\\

(i) There exists an invertible matrix $V$ with $\kap(V)\leq k$ such
that $D=V^{-1}AV$ is diagonal.

(ii) The functional calculus satisfies
\[
\norm f(A) \norm \leq k \norm f\norm_\infty
\]
for all complex-valued functions $f$ defined on $\Spec (A)$, where
\[
\norm f\norm_\infty = \max\{|f(\lam )|:\lam\in\Spec(A)\}.
\]
(iii) The resolvent operators satisfy
\[
\norm (A-z I)^{-1}\norm\leq k \,\,\dist (z,\Spec (A))^{-1}
\]
for all $z\notin \Spec (A)$, where $\dist$ is the Euclidean
distance of a point from a set.

(iv) The spectral projection $P$ of every eigenvalue $\lam$ of $A$
satisfies $\norm P\norm\leq k$, and hence
\[
i(A)\leq k.
\]

\end{proposition}

\begin{Proof}

(i) $\imp$ (ii) We use $f(A)=Vf(D)V^{-1}$, $\Spec (A)=\Spec(D)$
and
\[
\norm f(D) \norm = \norm f \norm_\infty.
\]

(ii) $\imp$ (iii) This is a matter of considering the particular
function $f(\lam)=(\lam-z)^{-1}$.

(iii) $\imp$ (iv) We express the spectral projection as a contour
integral of resolvent operators around a small circle centred at
$\lam$.
\end{Proof}

\begin{note}
The proof of the theorem remains valid if we replace the use of the
Euclidean norm on $\C^N$ by any other norm, provided the
appropriate operator norm is used for matrices and the operator
norm of any diagonal matrix $D$ is given by
\[
\norm D \norm =\max \{ |D_{n,n}|:1\leq n\leq N \}.
\]
This is equivalent to the norm being absolute \cite[sect. 10.5]{LT}
and holds in particular for all of the $l^p$ norms. However,
certain other matters, such as the identification of orthogonal
projections with those of norm $1$, are dependent on the use of the
Euclidean norm.
\end{note}

\begin{note}
If we assume condition (iv) of the above theorem then it follows
from the formula
\[
(A-z)^{-1}=\sum_i P_i(\lam_i-z)^{-1}
\]
that
\[
\norm (A-z)^{-1}\norm \leq Nk\,\,\dist(\Spec(A),z)^{-1}
\]
Since the value of $k$ frequently increases exponentially with the
dimension $N$,\cite{Tre1,triangular} one must expect pseudospectral
information and that obtained from the instability index to be
broadly equivalent.
\end{note}

The point of the theorem is that if any spectral projection of $A$
has very large norm, then the constant $k$ of any of the earlier
conditions must be very large, and diagonalization of the matrix
$A$ is an intrinsically ill-conditioned procedure.

There is a family of matrices for which the instability index
defined above can be computed in closed form. This is of some
interest for its own sake, but we used it to verify the algorithm
used to compute the instability indices of randomly generated
matrices. We assume that
\[
A_{m,n}=f(m-n)a^{m-n}
\]
where $a>1$ and $f:\Z\to \C$ is any function which is periodic with
period $N$.

\begin{theorem}
The eigenvectors of $A$ are of the form
\[
\phi_r(n)=a^n\rme^{2\pi i rn/N}
\]
where $r=1,...,N$. The corresponding spectral projections all have
norm
\[
c=\frac{a}{(a^2-1)N}\left( a^N-a^{-N}\right).
\]
Thus the instability index of $A$ is also $c$.
\end{theorem}

\begin{Proof} The first statement is a matter of applying the
matrix to such a vector, and noting that the set of all such
vectors is a basis for $\C^N$. The corresponding eigenvectors of
$A^\ast $ are
\[
\psi_r(n)=a^{-n}\rme^{2\pi i rn/N}.
\]
From these facts it is now easy to calculate the condition numbers
of each of the eigenvalues using (\ref{Pnorm}).
\end{Proof}

\section{Numerical Results}

We have applied the above ideas to a series of randomly generated
tridiagonal $N\times N$ matrices for various values of $N$. The
matrices are of the form $A_{m,n}$ so that $A_{m,n}=0$ for
$|m-n|>1$, the other coefficients being chosen randomly and
independently. If $m-n=1$ the coefficients are chosen using a
uniform distribution on $[0,1]$, if $m-n=0$ the coefficients are
chosen using a uniform distribution on $[0,2]$, and if $m-n=-1$ the
coefficients are chosen using a uniform distribution on $[0,3]$.

For each $A$ we used Matlab to obtain an invertible matrix $V$ and
a diagonal matrix $D$ such that $V^{-1}AV=D$. The columns of $V$
are then the eigenvectors $\phi_n$ of $A$, and are provided by
Matlab in normalized form. The rows of $V^{-1}$ are $\psi^\ast_m$
where $\psi_n$ are the eigenvectors of $A^\ast$. We get
$\psi^\ast_n\phi_n=1$ automatically, so
\[
\norm P_n \norm^2 = \norm \psi_n \norm^2 = E_{n,n}
\]
where $E=V^{-1}\left(V^{-1}\right)^\ast$. We thus finally get
\begin{equation}
i(A)=\max\{ |E_{n,n}|^{1/2}:1\leq n \leq N\}.\label{index}
\end{equation}

For each of $M$ randomly generated $N\times N$ matrices $A$ we
computed the instability index using (\ref{index}), and then sorted
the data points into increasing order. Defining $\cP_r$ to be the
number such that such that  $r\%$ of the instability indices were
less than $\cP_r$, we determined $\cP_{50}$ and $\cP_{95}$ for
various values of $N$. In fact we carried out each computation
twice in order to give some idea of the degree of reliability of
our results; we tabulated the average of the two values and the
difference $\cD_{r}$ expressed as a proportion of the average.

\[
\begin{array}{cccccccc}
N&M&\cP_{50}&\cD_{50}&\cP_{95}&\cD_{95}&\log(\cP_{50})/N&\log(\cP_{95})/N\\

\hline

10&10^6&18.47&0.007&304.1&0.001&0.292&0.572\\

20&10^6&684.9&0.005&6.13\times 10^4&0.008&0.326&0.551\\

30&4\times 10^5&2.138\times 10^4&0.013&7.48\times
10^6&0.001&0.332&0.528\\

40&10^5&5.77\times 10^5&0.047&6.68\times 10^8&0.050&0.332&0.508\\

50&10^5&1.43\times 10^7&0.016&4.60\times10^{10}&0.032&0.330&0.491

\end{array}
\]

It is clear from our results that both $\cP_{50}$ and $\cP_{95}$
increase extremely rapidly with $N$. In fact our results support
the conjecture that $\cP_{50}\sim\rme^{N/3}$ as $N\to\infty$.
Matlab is already having a little difficulty in computing the
eigenvalues of the matrices for $N=50$, and fails entirely for
$N=200$. Such exponential increase has also been found by Trefethen
et al for other models using pseudospectral
methods.\cite{Tre1,triangular}

An alternative approach to the above questions would be to compute
the expected value of $i(A)$ over a large sample of matrices $A$,
but this would have the disadvantage of being unduly influenced by
the very large size of $i(A)$ for a small proportion of choices of
$A$.

\section{Distribution of Norms of Spectral Projections}

Instead of studying $i(A)$ for randomly distributed matrices $A$,
one may examine how the norms of the individual spectral
projections are distributed. This may be done in several ways. In
the first we sort the norms of the spectral projections of a
particular $N\times N$ matrix $A$ in increasing order, but instead
of examining the largest of these, namely $i(A)$, we evaluate the
number $j(A)$ half way through the list. We did this for a series
of $10^5$ randomly generated $30\times 30$ matrices, and discovered
that for $50\%$ of these matrices $j(A)\leq 221.2$. The fact that
this number is so much smaller than $i(A)$ under the corresponding
conditions indicates that only a small proportion of the spectral
projections of a typical random matrix $A$ have really large norms.

A second procedure is to consider all of the norms of the spectral
projections of the $10^5$ randomly generated $30\times 30$ matrices
as one list of $30\times 10^5$ numbers. When we carried out this
numerical experiment we found that $50\%$ of the norms so obtained
were less than $157.8$. There is no reason why the above two
numbers should coincide, but they are of the same order of
magnitude, leading to the same conclusion.

The final and most interesting method is to carry out a spectral
analysis of the covariance matrix associated with the norms of the
spectral projections. We proceed as follows.

For each matrix $A$ we define the numbers $X_n\geq 0$ for $1\leq n
\leq N$ by
\begin{equation}
X_n=\log\left(\norm P_n\norm\right)   \label{X}
\end{equation}
where these are sorted in increasing order. (It is possible to
carry out similar calculations without taking the logarithm above,
but the results are less compelling.) As $A$ varies within the
usual class these provide a family of $N$ non-negative random
variables whose covariance matrix is defined by
\[
C_{m,n}={\bf E}\left[ X_mX_n\right].
\]
The eigenvalues and eigenvectors of this matrix provide information
about the distribution of the norms of the spectral projections of
`typical' random matrices.

We carried out the above computation for a sample of $10^5$
randomly distributed $10\times 10$ matrices. The eigenvalues of $C$
were found to be $0.0060$,
    $0.0196$,
    $0.0239$,
    $0.0323$,
    $0.0433$,
    $0.0613$,
    $0.1034$,
    $0.2133$,
    $0.7169$,
   $54.6547$.
The fact that one eigenvalue is so dominant is very striking, and
indicates that to a very good approximation most of the random
matrices have very similar distributions of the norms of their
spectral projections. We repeated the computation for a sample of
$10^5$ randomly distributed $30\times 30$ matrices. There were only
$4$ eigenvalues larger than $1$, these being $1.5$, $3.1$, $12.02$,
$1322.9$.

In both cases we computed the eigenvector $v$ corresponding to the
largest eigenvalue of the covariance matrix $C$. We found that
$v_n$ was close to being proportional to $n$. This corresponds to
the norms $\norm P_n
\norm$ of a `typical' random matrix $A$ increasing exponentially
with $n$. We conjecture that the dominance of the leading
eigenvalue increases and the rate of increase of the spectral norms
becomes more accurately exponential as the size $N$ of the matrices
increases.

Let us be more precise about this. Given $N$ we consider the sample
space
\[
\Ome_N=[0,1]^{N-1}\times [0,2]^N\times [0,3]^{N-1}
\]
provided with the uniform probability distribution. For each
$\ome\in\Ome_N$ we described how to construct a matrix $A_\ome$ and
then the $N$ random variables $X_n$ defined by (\ref{X}). These
have a symmetric $N\times N$ covariance matrix $C$ whose
eigenvalues may be ordered so as to satisfy
\[
0\leq \lam_1\leq...\leq\lam_N.
\]
We finally define
\[
\mu_N=\frac{\lam_{N-1}}{\lam_N}
\]
making explicit the dependence of $\mu$ on $N$. The conjecture is
then that
\[
\lim_{N\to\infty}\mu_N=0.
\]

We tested this hypothesis by considering a series of $T$ randomly
generated $N\times N$ matrices for various values of $T$ and $N$.
The results below provide some support for the conjecture.

\[
\begin{array}{ccccc}
N&T&\lam_{N-1}&\lam_N&\mu_N\\

\hline

10  &  10^6  &   0.725   & 54.99   & 0.0132  \\

20  &  10^6  &   4.422   &422.40   & 0.0105  \\

30  &  4\times 10^5  &   12.02   & 1319.5  & 0.0091  \\

40  &  10^5  &   24.15   & 2928.4  &  0.0082 \\

50  &  10^5  &   40.73   & 5408.9  & 0.0075

\end{array}
\]

\section{Operators with Randomly Distributed Coefficients}

In this section we present some ideas relating to a random family
of bounded linear operators acting on the infinite-dimensional
Hilbert space $l^2(\Z)$. This may be regarded as the infinite
volume limit of our earlier problems, although pseudospectral
theory suggests that one should also study the `same' operator on
$l^2(\Z^+)$. Physically the choice between these two operators
depends upon whether one wishes to include end effects, which are
present both for large finite intervals and for the operator on
$l^2(\Z^+)$. Our methods can easily be adapted to this case, but we
do not spell out the modifications needed.

We first formulate the ideas at a moderately general level, and
only later restrict attention to the non-self-adjoint Anderson
model of Hatano-Nelson. Let $A$ be the operator associated with an
infinite matrix $\{ A_{m,n}\}$ where $A_{m,n}=0$ if $|m-n|>1$; we
suppose that the vectors $v_n=\left(
A_{n+1,n},A_{n,n},A_{n,n+1}\right)\in\C^3$ are distributed
independently according to a common law $\mu$, where $\mu$ is a
probability measure on $\C^3$ with compact support $K$. It follows
from the assumptions and the fact that the matrix $\{ A_{m,n}\}$ is
tridiagonal that it is associated with a bounded operator $A$ such
that
\[
\norm A \norm \leq \sum_{r=1}^3 \max\{| w_r | :w\in K\}.
\]
The above procedure defines a stochastic family of bounded
operators $A_\ome$, for $\ome$ in the sample space $K^\infty$.  It
is of some interest that the results which we obtain are not truly
probabilistic: the statements of our theorems only involve the set
$K$ rather than the probability measure $\mu$.

In order to prove some results about the spectra of operators in
the family we introduce a notion from \cite{DS}. Given any bounded
operator $X$ on a Hilbert space $\cH$, we say that the operator $Y$
lies in its limit class, $Y\in\cC(X)$, if there exists a sequence
$U_s$ of unitary operators on $\cH$ such that $U_s^\ast X U_s$
converges strongly to $Y$ as $s\to\infty$. We also define the
approximate point spectrum $\sig(X)$ of $X$ to be the set of all
$\lam\in\C$ for which there exists a sequence of vectors
$f_s\in\cH$ such that $\norm f_s \norm =1$ and
\[
\lim_{s\to\infty}\norm Xf_s-\lam f_s\norm =0.
\]
In all the examples in this paper $\sig(X)=\Spec(X)$ almost surely
but it seems desirable for the sake of possible future applications
to keep the logical distinction. The following two lemmas are
modifications of a classical result of Pastur stating that the
spectrum of random tridiagonal operators is almost surely constant;
the proof uses translation ergodicity of the class of
operators.\cite[p 167]{Cycon}

\begin{lemma} If $Y$ lies in the limit class of $X$ then
$\sig(X)\supseteq \sig (Y)$. Hence
\[
\Spec(X)\supseteq\bigcup\{ \sig(Y):Y\in\cC(X)\}.
\]
In particular if each lies in the limit class of the other then
$\sig(X)=\sig(Y)$.
\end{lemma}

\Proof Suppose that for some $\lam \in\C$ and all $\eps >0$
there exists $f\in\cH$ such that $\norm f \norm =1$ and $\norm Y
f-\lam f\norm <\eps$. Now put $f_s=U_sf$ and observe that
\[
\norm Xf_s-\lam f_s \norm =\norm U_s^\ast X U_s f-\lam f\norm\to
\norm Y f-\lam f\norm
\]
as $s\to\infty$. Therefore $\norm Xf_s-\lam f_s \norm <\eps$ for
all large enough $s$, and $\lam\in\sig (X)$.

We apply the above to the randomly generated operators $A_\ome$
acting on $l^2(\Z)$.

\begin{lemma}
The limit class of almost every operator $A_{\ome}$ generated as
described contains every operator $A_{\tilde\ome}$ such that
$\tilde\ome\in K^\infty$.
\end{lemma}

\Proof Let $N\in \Z^+$ and $\tilde\ome =\{ \tilde v_n\}_{n\in\Z} \in
K^\infty$. Given $\eps >0$ we put
\[
V_n=\{ z\in\C^3:\norm \tilde v_n-z\norm <\eps\}
\]
so that $\mu(V_n)>0$ for all $n$. Taking $U_s:l^2(\Z)\to l^2(Z)$ to
be the unitary operators associated with appropriate translations
of $\Z$, we need to show that for almost every $\ome=\{v_n \} \in
K^\infty$ there exists $M$ such that $v_{n+M}\in V_n$ for all
$-N\leq n\leq N$.

To prove this we put $V=\prod_{n=-N}^N V_n$ and
\[
\ome=\{v_n\}_{n\in\Z}=\{w_m\}_{m\in\Z}
\]
where
\[
w_m=\{ v_{m(2N+1)+r}\}_{r=-N}^N\in K^{2N+1}.
\]
The vectors $w_m$ are independent and identically distributed with
positive probability that $w_m\in V$. Hence the probability that
none of the $w_m$ lie in $V$ is zero.

As an application of the lemma, let $v,w\in C^3$ and let $B_{v,w}$
be the bounded operator associated with the infinite matrix $\{
B_{m,n}\}$ such that $B_{m,n}=0$ if $|m-n|>1$; we also assume that
the vector $\left( B_{n+1,n},B_{n,n},B_{n,n+1}\right)$ equals $v$
if $n\geq 1$ and equals $w$ if $n\leq 0$. The point of introducing
this class of operators is that they are similar to operators whose
spectrum is well understood.\cite{Boe,Boe2,Kre}

\begin{corollary}
With probability $1$ one has
\[
\Spec (A_\ome )\supseteq \bigcup \left\{ \sig(B_{v,w}):v,w\in K\right\}.
\]
\end{corollary}

\Proof Every operator $B_{v,w}$ lies in the limit class of $A$ by
Lemma 5.

The significance of the above lemmas is best seen by applying them
to an example. We assume that $\mu$ is a probability measure
concentrated on a finite subset $F$ of $\R^2$ (more general
probability measures can also be treated). We assume that
$\mu(0,0)>0$ and that any other point $(x,y)$ with $\mu(x,y)>0$
satisfies $x>0$ and $y>0$. We then define the random family of
operators $A$ on $l^2(\Z)$ as described above with the following
simplification: we assume that all $A_{n,n}=0$ and that the vectors
$\left( A_{n+1,n},A_{n,n+1}\right)$ are distributed independently
according to the law of $\mu$.

\begin{theorem} Depending to the choice of $\mu$, either $A_\ome$
is self-adjoint with probability $1$ or it is non-self-adjoint with
probability $1$. In the latter case $A_\ome$ possesses a countable
set of eigenvalues whose corresponding eigenvectors span a dense
linear subspace of $l^2(\Z)$. With probability $1$ the set of
eigenvectors is not a basis of $l^2(\Z)$.
\end{theorem}

\Proof  The self-adjoint case occurs when the support of $\mu$ is
contained in the diagonal set $\{ (x,y)\in\R^2:x=y\}$ and we assume
that this is not the case below.

With probability $1$ an infinite number of the pairs
$(A_{n+1,n},A_{n,n+1})$ are equal to $(0,0)$, and we assume that
this happens for the increasing sequence $\{ N_r\} $ where
$r\in\Z$. It may then be seen that $A_\ome$ can be decomposed as
the orthogonal direct sum of matrices $C_r$ of sizes $M_r\times
M_r$ where $M_r=N_r-N_{r-1}$. Since each matrix $C_r$ is
tridiagonal with positive off-diagonal entries and zero diagonal
entries, its eigenvalues are all real and of multiplicity one. By
combining all of the eigenvectors of the $C_r$ as $r$ increases we
see that the eigenvectors of $A_\ome$ almost surely span a dense
linear subspace of $l^2(\Z)$. It remains only to prove that this
set of eigenvectors is almost surely not a basis.

Let $(x,y)$ be any point in $\R^2$ with $0<x<y$ and $\mu(x,y)>0$.
Since the pairs of coefficients of $A_\ome$ are chosen
independently, among the $C_r$ there almost surely exist all
$s\times s$ matrices of the form $B_s$ where
\[
B_s=\left\{ \begin{array}{ll}

x&\mbox{if $\,\,m=n+1$}\\ y&\mbox{if $\,\,m=n-1$}\\
0&\mbox{otherwise}.

\end{array} \right.
\]
Now the eigenvectors of each $B_s$ are given explicitly by
\[
\phi_{k,r}=\sin\left( \frac{\pi kr}{s+1}\right)
\left( \frac{x}{y}\right)^{r/2}
\]
and the corresponding eigenvectors of the adjoint operator are
\[
\psi_{k,r}=\sin\left( \frac{\pi kr}{s+1}\right)
\left( \frac{y}{x}\right)^{r/2}
\]
where $k$ labels which eigenvector is being considered and $r$
which coefficient.

It follows that
\begin{eqnarray*}
\langle \phi_k,\psi_k\rangle&=&\sum_{r=1}^s\sin\left( \frac{\pi kr}{s+1}\right)^2\\
\norm \phi_k \norm^2&=&\sum_{r=1}^s\sin\left( \frac{\pi kr}{s+1}\right)^2\left|\frac{x}{y}\right|^r\\
\norm \psi_k \norm^2&=&\sum_{r=1}^s\sin\left( \frac{\pi kr}{s+1}\right)^2\left|\frac{y}{x}\right|^r
\end{eqnarray*}

An application of (\ref{Pnorm}) now shows that
$\lim_{s\to\infty}i(B_s)=+\infty$, from which it follows that with
probability $1$ the norms of the spectral projections of $A_\ome$
are not uniformly bounded. This implies that the eigenvectors
cannot form a basis.

An example of Zabzyk which has some similarities to those of the
above theorem is discussed in Theorem 2.17 of \cite{OPS}, where the
failure of the spectral mapping theorem for semigroups is
demonstrated. Another type of example involving differential
operators whose eigenvectors do not form a basis was presented in
\cite{Dav2,Dav4,Dav3}. It appears that such a situation is
relatively common for non-self-adjoint operators in infinite
dimensions. Since the eigenvectors do not form a basis for the
Hilbert space, there is no reason to expect that they determine the
spectral behaviour of $A_\ome$. We continue with the hypotheses
formulated before Theorem 7.

\begin{theorem}
The eigenvalues of the operator $A_\ome$ are almost surely all
real. However the spectrum of $A_\ome$ almost surely contains the
interior of the ellipse $E$ defined by
\[
E=\{ x\rme^{i\theta}+y\rme^{-i\theta}:\theta\in [0,2\pi]\}
\]
for every non-zero $(x,y)\in F$.
\end{theorem}

\Proof
The operator $A_\ome$ is almost surely the orthogonal direct sum of
tridiagonal matrices $C_r$, where we continue to use our previous
notation. Each $C_r$ is similar to a real symmetric matrix and
therefore has real eigenvalues.

If $v=(x,y)$ is a non-zero point in $F$ then by Corollary 6 the
spectrum of $A_\ome$ almost surely contains $\sig(B_{0,v})$. We
assume that $0<x<y$, the case $0<y<x$ having a similar analysis
involving the adjoint operators. A direct computation shows that if
$z$ lies inside the said ellipse then  both solutions $w_i$ of
\[
x-zw+yw^2=0
\]
satisfy $|w_i|<1$. It follows that $z$ is an eigenvalue of
$B_{0,v}$, the eigenvector being the sequence $f\in l^2(\Z)$ given
by
\[
f_r=\left\{ \begin{array} {cl} w_1^r-w_2^r&\mbox{ if $r\geq 1$}\\
0&\mbox{otherwise}.
\end{array}
\right.
\]

\section{The Non-Self-Adjoint Anderson Model}

In this section we consider non-self-adjoint Anderson-type
operators of the form
\begin{equation}
Hf_n=\rme^{-g}f_{n-1}+\rme^gf_{n+1}+V_nf_n\label{schr}
\end{equation}
where $g>0$ and $V$ is a random real-valued potential. We assume
that the values of $V$ at different points are independent and
identically distributed according to a probability law $\mu$ which
has compact support $M\subseteq \R$.

A considerable amount of attention has already been paid to such
operators, which arise in population biology and solid state
physics, as cited in \cite{Gold,Hatano1,Hatano2,Nelson}. If one
supposes that the operator acts on $l^2\{-N,N\}$ subject to
Dirichlet boundary conditions then its spectrum is almost surely
real, because $H$ is similar to the operator defined by the same
formula but with $g=0$. The similarity is determined by the
operator
\[
Sf_n=\rme^{ng}f_n
\]
which is bounded on this finite-dimensional space. On the other
hand if one imposes periodic boundary conditions, then the spectrum
of $H$ is much more interesting and has been analyzed in great
detail in the limit of large $N$ both numerically and
theoretically.\cite{Gold,Hatano1,Hatano2,Nelson} The imposition of
periodic boundary conditions was justified in Nelson-Shnerb by the
fact that they were considering a rotationally invariant problem
(the diffusion of bacteria in a rotating nutrient). Goldsheid has
observed that the spectral results proved in \cite{Gold} for
periodic boundary conditions are equally valid for a wide range of
quasi-periodic boundary conditions. However, Trefethen has pointed
out that the spectral behaviour of this type of operator is highly
problematical.\cite{Trefcomm} This phenomenon has been investigated
from several points of view over the last decade,\cite{Boe, Boe2,
Dav1,Red, RT, Tre1,Tre2} and among the conclusions is the warning
that one cannot assume that a solution of a non-linear equation is
stable simply because the eigenvalues of its linearization about
the solution have negative real parts.

We show that if one considers the non-self-adjoint Anderson model
$H$ acting on $l^2(\Z)$, the spectrum is entirely different from
what one obtains by letting $N\to\infty$ subject to any of the
above types of boundary condition. The following theorems are
directed towards locating the spectrum of $H$, but it would clearly
be highly desirable to find a precise formula for it.

\begin{theorem}\label{inside}
Let $H$ be defined on $l^2(\Z)$ by (\ref{schr}) where $V$ satisfies
the stated conditions. Then
\[
\Spec(H)\subseteq  \conv(E)+\conv(M)
\]
almost surely, where $\conv$ denotes convex hull and $E$ is the
ellipse
\begin{equation}
E=\{\rme^{g+i\theta}+\rme^{-g-i\theta}:\theta\in[0,2\pi]\}.\label{ellipse}
\end{equation}
\end{theorem}

\Proof Let $A$ be the operator obtained from $H$ by deleting the
potential $V$, so that $\Spec(A)=E$. Then the result follows from
the equation
\[
\Spec(H)\subseteq \Num(H)\subseteq
\Num(A)+\Num(V)=\conv(E)+\conv(M).
\]
where $\Num$ stands for the numerical range
\[
\Num (A)=\{ \langle Af,f\rangle:\norm f \norm =1\}.
\]

\begin{theorem}\label{pert}
Under the above assumptions we also have
\[
\Spec(H)\subseteq \{z\in\C:\dist(z,E)\leq m\}
\]
almost surely, where $\dist$ is the distance function and
\[
m=\max\{|s|:s\in M\}.
\]
\end{theorem}

\Proof Since $A$ is a normal operator with spectrum equal to $E$, we have
\[
\norm (z-A)^{-1}\norm =\dist(z,E)^{-1}
\]
for all $z\notin E$. Also $\norm V\norm =m$ almost surely, and the
result follows by examining the convergence of the perturbation
series
\[
(z-A-V)^{-1}=\sum_{n=0}^\infty
(z-A)^{-1}\left(V(z-A)^{-1}\right)^n.
\]

In the reverse direction we may apply Lemma 4 and Theorem 5 to
obtain the following result, which is further improved in Theorem
\ref{sigl}.

\begin{theorem}\label{outside}
With probability one, the spectrum of $H$ contains the set $E+m$
for every $m\in M$.
\end{theorem}

\Proof The operator $A+mI$ lies in the limit class of $H$ for all
$m\in M$ and its approximate point spectrum is $E+m$ (as is its
spectrum).

\begin{example}
Let $\mu$, the measure determining $V$, have support $[-B,B]$. By
applying Theorems \ref{inside} and \ref{outside} it follows that
the spectrum of $H$ is almost surely equal to the convex set
$E+[-B,B]$ provided $B\geq
\rme^g+\rme^{-g}$.  This is quite different from what occurs for
the same problem on $l^2(-N,N)$ in the limit $N\to\infty$, whatever
boundary conditions are assumed. The point is that for such
operators approximate eigenvalues need not be close to genuine
eigenvalues, and the full spectral behaviour of the operators can
best be seen using pseudospectral theory.

The determination of the spectrum for smaller values of $B$ is more
complicated. An application of Theorem \ref{pert} shows that if
$B<\rme^g-\rme^{-g}$ then
\[
\Spec(H)\cap \{z\in\C:|z|<r\}=\emptyset
\]
where $r=\rme^g-\rme^{-g}-B>0$. Since $E\subseteq \Spec(H)$, it
follows that the spectrum of $H$ has a hole in it.

We have not been able to determine the precise range of values of
$g,B$ for which the spectrum contains a hole. A numerical solution
seems out of the question because the spectrum is heavily affected
by extremely infrequent `regular structure' in the potentials.
\end{example}

\section{Classification of the Spectrum}

The classification of the spectrum of non-self-adjoint operators is
in a primitive state by comparison with that of self-adjoint
operators. We start at an abstract level.

If $f\in l^2(\Z)$ satisfies $\norm f \norm =1$, we put
\[
\var (f) = \langle Q^2 f,f\rangle -\langle Qf,f\rangle^2
\]
where $Q$ is the position operator, provided this is finite. If $A$
is a bounded linear operator acting on $l^2(\Z)$, we say that
$\lam$ lies in its localized spectrum $\sigl(A)$ if there exists
$c$ and a sequence $f_n\in l^2(\Z)$ such that $\norm f_n
\norm =1$ and $\var(f_n)\leq c$ for all $n$ and
\[
\lim_{n\to\infty} \norm Af_n-\lam f_n\norm =0.
\]
A more general definition is possible but not necessary for our
purposes. The localized spectrum of $A$ includes all eigenvalues
whose associated eigenvectors have finite variances, but need not
be a closed set. We say that $B$ lies in the translation limit
class of $A$ if there exists a sequence of unitary translation
operators $U_n$ such that $U_n^\ast AU_n$ converges strongly to $B$
as $n\to\infty$.

\begin{lemma}
If $B$ lies in the translation limit class of $A$ then
$\sigl(A)\supseteq \sigl(B)$. In particular $\lam\in\sigl(A)$ for
every eigenvalue $\lam$ of $B$ whose corresponding eigenvector has
finite variance.
\end{lemma}

\Proof This is the same as that of Lemma 5, with the extra
observation that the variance is unchanged by translations.

The following theorem establishes that the localized spectrum is
distinct from the approximate point spectrum. We say that an
operator $A$ acting on $l^2$ is a finite order convolution operator
if it is of the form $Af=k\ast f$ where $\ast$ denotes convolution
and $k$ has finite support.

\begin{theorem}
If $H=A+B$ where $A$ is a finite order convolution operator and $B$
is compact, then its localized spectrum is the set of all its
eigenvalues whose corresponding eigenvectors have finite variance.
If $B=0$ this set is empty.
\end{theorem}

\Proof Let $\lam\in\sigl(H)$ and let $f_n$ satisfy $\norm f_n\norm
=1$, $\var(f_n)\leq c$ and $\norm Hf_n-\lam f_n\norm\to 0$ as
$n\to\infty$. If we put $a_n=\langle Qf_n,f_n\rangle$ then we have
two cases to consider.

If $a_n$ is an unbounded sequence then by passing to a subsequence
we may assume that $a_n$ diverges. If we define $U_n$ by
\[
(U_n\phi)_m=\phi_{m-a_n}
\]
and put $g_n=U_n^\ast f_n$  then $\langle Qg_n,g_n\rangle=0$ and
$\var(g_n)=\var(f_n)$. Hence $g_n$ lies in the compact set
\[
S=\{ \phi \in l^2(\Z):\norm \phi\norm =1 {\rm \,\, and\,\,} \langle
Q^2 \phi,\phi\rangle
\leq c\}
\]
for all $n$. By passing to a convergent subsequence we assume
further that $g_n$ converges to a limit $g\in S$. Now we also have
\[
\lim_{n\to\infty} \norm U_n^\ast H U_n g_n -\lam g_n\norm =0
\]
where $U_n^\ast H U_n$ converges strongly to $A$ as $n\to\infty$.
Hence $A g=\lam g$, which contradicts the fact, proved using
Fourier analysis, that the point spectrum of $A$ is empty.

The alternative case is that there exists $a$ such that $|a_n|\leq
a$ for all $n$. It then follows that $f_n$ lies in the compact set
\[
T=\{ \phi \in l^2(\Z):\norm \phi\norm =1 {\rm \,\, and\,\,} \langle
Q^2 \phi,\phi\rangle
\leq c+a^2\}
\]
for all $n$. By passing to a convergent subsequence we assume
further that $f_n$ converges to a limit $f\in T$. It is now
immediate that $Hf=\lam f$.

The final special case of the theorem follows from the fact that
$A$ has empty point spectrum.

In spite of the above, in the following context the localized
spectrum is quite different from the point spectrum. We place
ourselves in the situation described in the first paragraph of
Section 6. Using the formula (\ref{ellipse}) for the ellipse $E$,
we define $I_z$ to be the open interior of $E+z$ and $O_z$ to be
the open exterior of $E+z$ for any $z\in\C$.

\begin{theorem}\label{sigl}
If $H$ is defined by (\ref{schr}) then
\[
\sigl(H)\supseteq \bigcup_{\alp,\bet\in M}\{
I_{\alp}\cap O_{\bet}\}
\]
almost surely.
\end{theorem}

\Proof Given $\alp,\bet\in M$ the operator $K$ defined by
\[
(Kf)_n=\left\{ \begin{array} {cl}
\rme^{-g}f_{n-1} +\alp f_n +\rme^g f_{n+1}&\mbox{if $n\geq 1$}\\
\rme^{-g}f_{n-1} +\bet f_n +\rme^g f_{n+1}&\mbox{if $n\leq 0$}
\end{array}\right.
\]
lies in the translation limit class of $H$ by Lemma 5, so it is
sufficient to prove that if $\lam \in I_{\alp}\cap O_{\bet}$ then
$\lam$ is an eigenvalue of $K$ whose corresponding eigenvector has
finite variance.

The eigenvalues and eigenvectors of $K$ may be determined
explicitly. For $n\geq 2$ we solve the equation
\[
\rme^{-g}u^{-1} +\alp +\rme^g u=\lam
\]
and observe using the Residue Theorem that $\lam\in I_\alp$ implies
that there are two solutions $u_1,u_2$ which both satisfy
$|u_i|<1$. For $n<0$ we solve
\[
\rme^{-g}v^{-1} +\bet +\rme^g v=\lam
\]
and observe similarly that $\lam\in O_\bet$ implies that there is a
solution $v$ which satisfies $|v|>1$. We now define
\[
f_n=\left\{ \begin{array} {cl} v^n&\mbox{if $\,\,n\leq 0$}\\
c_1u_1^n+c_2u_2^n&\mbox{if $\,\,n\leq 1$}.
\end{array}\right.
\]
By an appropriate choice of $c_1,c_2$ we can ensure that $Kf=\lam
f$ as required.

\section{Higher Dimensions}

The results of the last sections can all be extended to higher
dimensions, and we briefly describe the situation in $l^2(\Z^2)$.
We assume that
\[
Hf_{m,n}=\rme^{-g}f_{m-1,n}+\rme^{g}f_{m+1,n}
+\rme^{-h}f_{m,n-1}+\rme^{h}f_{m,n+1}+V_{m,n}f_{m,n}
\]
where $g>0$, $h>0$ and $V$ is a random real-valued potential. We
assume that $V$ has independent values at different points and that
they are distributed according to a probability measure $\mu$ with
compact support $M\subseteq \R$. If $V=0$ then Fourier analysis
shows that
\[
\Spec(H)=E+F
\]
where $E,F$ are the ellipses
\[
E=\{  \rme^{g+i\theta}+\rme^{-g-i\theta}:\theta\in [0,2\pi] \}
\]
and
\[
F=\{  \rme^{h+i\theta}+\rme^{-h-i\theta}:\theta\in [0,2\pi] \}.
\]

It is routine to show that Theorem \ref{outside} has the analogue

\begin{theorem}
With probability one, the spectrum of $H$ contains the set $E+F+m$
for every $m\in M$.
\end{theorem}

Localization of the spectrum is more complicated but the following
ideas provide some information. Let $A,B$ be two bounded
real-valued functions defined on $\Z$ and define
\[
(A+B)_{m,n}=A_m+B_n.
\]
\begin{theorem} If $A+B$ lies in the translation limit class of
$V$ almost surely, then
\[
\sigl(H)\supseteq \sigl(H_1)+\sigl(H_2)
\]
almost surely, where $H_j$ act on $l^2(\Z)$ according to
\[
H_1f_n=\rme^{-g}f_{n-1}+\rme^gf_{n+1}+A_nf_n
\]
and
\[
H_2f_n=\rme^{-h}f_{n-1}+\rme^hf_{n+1}+B_nf_n.
\]
\end{theorem}

\Proof We first observe that $H_1\tensor I+I\tensor H_2$ lies in
the translation limit class of $H$ almost surely, so
\[
\sigl(H)\supseteq \sigl(H_1\tensor I+I\tensor H_2).
\]
By considering test functions of the form $f_1\tensor f_2$ we also
see that
\[
\sigl(H_1\tensor I+I\tensor H_2)\supseteq \sigl(H_1)+\sigl(H_2).
\]

It might seem that the hypothesis of this theorem is rather
special. However if the support $M$ of $\mu$ is an interval then
the theorem is of real value. Under these assumptions if
$\alp,\bet\in M $ then $A+B$ lies in the translation limit class of
$V$ almost surely, where
\[
A_m=B_m=\left\{\begin{array}{cl}
\alp/2 &\mbox{if $m\geq 0$}\\
\bet/2 &\mbox{otherwise.}
\end{array}
\right.
\]
The localized spectrum of such operators is considered in Theorem
\ref{sigl}.

\section{Conclusions}

We have analyzed the spectral behaviour of a family of randomly
generated non-self-adjoint matrices by a variety of different
methods. Our conclusion is that their eigenvalues depend very
sensitively on the matrix entries even for quite small matrix
sizes. The standard proofs of the existence of the eigenvalues
depend upon the fundamental theorem of algebra, of which there are
many proofs, but it is clear that the well-known instability of the
roots of high degree polynomials leads to the large value of the
instability index of general non-self-adjoint matrices. This
problem does not occur for self-adjoint matrices because the
variational theorem implies that the eigenvalues of such matrices
do not change much under small perturbations, \cite{STDO}. Our
results indicate that nothing of a comparable nature is likely to
be available in the non-self-adjoint case.

We have also investigated a family of randomly generated
non-self-adjoint bounded operators acting on an
infinite-dimensional Hilbert space, for which the eigenvectors
almost surely generate a dense linear subspace. In spite of this
the eigenvectors almost surely do not form a basis and the
eigenvalues almost surely generate only a small part of the
spectrum. We finally proved that the full spectrum of the infinite
volume non-self-adjoint Anderson model bears little relationship
with the infinite volume limit of the spectra of the same operator
in finite intervals.

\vskip 1in

%%%%%%%%%%%%%%%%%%%%%%%%%
{\bf Acknowledgments } I would like to thank L N Trefethen, M
Embree and I Goldsheid for valuable comments.
%\par
%\vskip 1.5in
%%%%%%%%%%%%%%%%%%%%%%%

%%%%%%%%%%%%%%%%%%%%%%%%%%%%%%%%%%%%%%%
%\newpage
%%%%%%%%%%%%%%%%%%%%%%%%%%%%%%%%%%%%%%%%
Department of Mathematics\newline%
King's College\newline%
Strand\newline%
London WC2R 2LS\newline%
England
\vskip 0.2in
e-mail: E.Brian.Davies@kcl.ac.uk
\vfil
\end{document}